\documentstyle[12pt]{article}

\input epsf
\def\frak#1{{\bf #1}}
\def\Bbb#1{\bf#1}
\def\bb#1{\bf#1}

\def\llcorner#1{(#1)}
\newtheorem{example}{Example}[subsection]

\date{}

\newtheorem{Definition}{Definition}[subsection]
\newtheorem{Proposition}{Proposition}[subsection]
\newtheorem{proposition}{Proposition}[subsection]

\def\textmap#1{\mathop{\vbox{\ialign{
				##\crcr
    ${\scriptstyle\hfil\;\;#1\;\;\hfil}$\crcr
    \noalign{\kern-1pt\nointerlineskip}
    \rightarrowfill\crcr}}\;}}

\begin{document}

\begin{center}
\vspace{0.7 cm}
{\large{\bf Groups of Flagged Homotopies and Higher Gauge Theory}
\footnote{This work was partially supported by INTAS-OPEN-97-2072 grant.}} \\
\vspace{1.5 cm}
 Valery DOLOTIN \\
\vspace{0.2 cm}
{\it vd@mccme.ru } \\
\vspace{1.4 cm}
\end{center}

\begin{abstract}
Groups $\Pi_k(X;\sigma)$ of "flagged homotopies" are introduced of which the 
usual (Abelian for $k>1$) homotopy groups $\pi_k(X;p)$ is the limit case for 
flags $\sigma$ contracted to a point $p$. The calculus of exterior forms
with values in an algebra $A$ is developped of which the limit cases are 
the differential forms calculus (for $A=\bb R$) and gauge theory (for 1-forms). 
Moduli space of integrable forms with respect to higher gauge transforms 
(cohomology with coefficients in $A$) is introduced
with elements giving representations of $\Pi_k$ in $G=\exp A$.
\end{abstract}

\section{Basic Example}

Let $U$ be a local chart of a real 2-dimensional manifold $X$. 
Take a
nondegenerate map $\varphi :U\rightarrow {\Bbb C}^1$ providing $U$ with a
complex structure. For each map $f:U\to \bb C$ which is holoorphic with
respect to the complex structure above we have the map: 
$\exp f:U\rightarrow {\Bbb C}^{*}$. 
The differential of $\exp f$ is a map 
$A:TU\rightarrow T{\Bbb C}^{*}$, where $T{\Bbb C}^{*}=:g({\Bbb C}^*)=\Bbb C$
\negthinspace is the Lie algebra of the group ${\Bbb C}^{*}$. 
In the local chart $U\times 
{\frak g}({\Bbb C}^{*})$ of rank 2 vector bundle the \negthinspace 
$d(\exp f)$ may be written as the connection form $dy_i+\sum%
\limits_{j,k}A_{ij}^ky_jdx_k$ with coefficients \negthinspace $A^1(x_1,x_2)=\left( 
\begin{array}{cc}
u & v \\ 
-v & u
\end{array}
\right) \ {\rm and}\ 
A^2(x_1,x_2)=\left( 
\begin{array}{cc}
v & -u \\ 
u & v
\end{array}
\right) $, where $(x_1,x_2)$ are real and imaginary parts of 
complex structure coordinates
induced by $\varphi $, and $(y_1,y_2)$ \negthinspace are coordinates in the
fiber ${\frak g}({\Bbb C}^{*})$. This connection is automaticaly flat, since
the zero-curvature conditions $\frac{\partial A_i}{\partial x_j}-\frac{%
\partial A_j}{\partial x_i}+[A_i,A_j]$\negthinspace $=0$ in this case are
just the Cauchy-Riemann conditions \negthinspace $\frac{\partial u}{\partial x_2}%
=\frac{\partial v}{\partial x_1},\frac{\partial u}{\partial x_1}=-\frac{%
\partial v}{\partial x_2}$. So for a given complex structure we get a family
of flat connections corresponding to holomorphic maps.
Now vice versa, if we have a flat connection $A$
on a trivial real rank 2 vector bundle over a real 2-dimensional manifold $X$,
then, for each local chart $U\times F$ and fixed base point $x'\in X$, 
we may integrate A along paths $\gamma$ with $\gamma(0)=x'$, wich gives us a map
$U\to GL(2,{\Bbb R}),\ x\mapsto G_x:=\int\limits_{\gamma(x',x)}e^A$. If we 
fix fiber point $y'\in F_{x'}$ then we get a map 
\negthinspace $\psi :U\rightarrow F\simeq {\Bbb R}^2, 
x\mapsto G_x(y')$. In particular $\psi(x')=y'$.
Now, if we fix a complex structure on $F$ and $A$ is nondegenerate on $U$,
then $\psi $ induces a complex structure on $U$. A change of coordinates on $%
F$ which preserves complex structure may be represented as an action of
the operator \negthinspace $\left( 
\begin{array}{cc}
a & b \\ 
-b & a
\end{array}
\right) $. Then a global smooth change of fiber coordinates preserving
complex structure on each \negthinspace $F_x$ is the smooth function $\varphi
:x\mapsto \left( 
\begin{array}{cc}
a(x) & b(x) \\ 
-b(x) & a(x)
\end{array}
\right) $, 
giving a ''gauge field'' \negthinspace $g(x)=\exp \varphi
:X\rightarrow {\Bbb C}^{*}$. Such a change of fiber coordinates induces the
transform of connection coefficients \negthinspace $g(A)=A+dg\cdot g^{-1}$.
So given a connection \negthinspace $A_0$ we get a map of ''gauge group'' $%
{\cal G}$ to the space ${\cal A}$ of connections, such that \negthinspace $%
{\rm Id}\mapsto A_0$ and ${\cal A}$ gets fibered into orbits of ${\cal G}$-action. 

{\it Question:} when a flat connection $A$ may be obtained as the differential of a
globally defined map \negthinspace $X\rightarrow F={\Bbb C}^1$, or what is
the space ${\cal A}^{\prime }/{\cal G}$ of ${\cal G}$-fibers in the flat
connections subspace ${\cal A}^{\prime }\subset {\cal A}$? In particular,
when $X$ is orientable and has a fixed complex structure (Riemann surface),
then it may be repersented by various connections, which, in general, belong to
different ${\cal G}$-fibers. The integration of $A$ with given $v\in x\times
F$ in this case is called analytic continuation. 


\section{General construction}

\subsection{Space of Flagged Homotopies}

{\it Step 1}

Take a $1$-dimensional path $\gamma ^1:[0,1]\rightarrow {\Bbb R}^n$. Its
homotopy $\gamma ^2$ with the fixed boundary $\sigma _0=\gamma ^1(0),\sigma
_1=\gamma ^1(1)$ sweeps a $2$-disk $D$. The boundary $\partial D=S^1$ gets a
decomposition into pairs of cells $((\sigma _0,\sigma _1),(\sigma
_0^1,\sigma _1^1))$, where $(\sigma _0,\sigma _1)$ are the boundaries of $\gamma
^1$ (points) and $(\sigma _0^1,\sigma _1^1)$ are the boundaries of the homotopy $%
\gamma ^2$ (semi-circles forming $S^1=\partial D$) (see Fig.1). This
decomposition will be called {\em flagging} of $D$. For each point $p\in D$
we have two numbers $t_1(p),t_2(p)$ - the values at $p$ of parameters of
homotopies $\gamma ^1$ and $\gamma ^2$ correspondingly. Thus we have a
coordinatization of $D$.

\begin{center}
\epsfxsize=15cm

Fig.1
\end{center}

\noindent
{\it Inductive step}

Take a $k$-dimensional disk $D^k$ and take its homotopy $\gamma ^{k+1}$ with
the fixed boundary $\partial D^k=S^{k-1}$. It sweeps a $k+1$-disk $D$. If $%
D^k$ is flagged then $\partial D=S^k$ gets flagging with $k$-cells $(\sigma
_0^k,\sigma _1^k)$ being boundaries of $\gamma ^{k+1}$. For each point $p\in
D$ we have the value $t_{k+1}(p)$ of parameter of homotopy $\gamma ^{k+1}$ such
that $p\in \gamma ^{k+1}(t(p))$. If we fix the homotopy on $\sigma _0^k$,
then the corresponding $t^k$-parameter function on $\sigma _0^k$ is extended
via homotopy $\gamma ^{k+1}$ to the $t^k$-coordinatization on $D$. So fixing
a flag of homotopies $(\gamma ^1,\gamma ^2,\ldots ,\gamma ^k)$ on $(\sigma
_0^1,\ldots ,\sigma _0^k)$ correspondingly, where the body $\sigma _0^i$ of
each $\gamma ^i=\gamma ^{i+1}(0)$ is the initial point of the subsequent
flag element $\gamma ^{i+1}$, gives us values $(t^1(p),\ldots ,t^{k+1}(p))$
of the parameters of these homotopies for each $p\in D$ and thus
coordinatization on $D$.

For a manifold $X$ the set of homotopies $\gamma ^k$ with fixed flagged
boundaries will be denoted by $\Gamma ^k(X)$. In particular, $\Gamma ^0(X)=X$%
. For a flagging $\sigma =((\sigma _0,\sigma _1),\ldots ,(\sigma _0^d,\sigma
_1^d))$ the set of elements of $\Gamma ^k(X)$ having this flagging is
denoted by $\Gamma _\sigma ^k(X)$. Elements of $\Gamma _\sigma ^k(X)$ are
subsets of $X$ on one hand and paths in $\Gamma _\sigma ^{k-1}(X)$ on the
other. Then for $\gamma ^k\in \Gamma _\sigma ^k(X),\gamma \in \Gamma _\sigma
^{k+1}(X)$ we may write $\gamma ^k\subset \gamma $ or $\gamma ^k\in \gamma $
correspondingly.

\subsubsection{Composition Law}

Take a pair $\gamma ^{\prime },\gamma ^{\prime \prime }\in \Gamma _\sigma ^k$
such that $\sigma _0^k(\gamma ^{\prime \prime })=\sigma _1^k(\gamma ^{\prime
})$. Then we get an element $\gamma ^{\prime \prime }\circ \gamma ^{\prime
}\in \Gamma _\sigma ^k$. This composition law makes $\Gamma _\sigma ^k$ into
a category called {\em space of flagged homotopies}.

\subsubsection{Flagged Homotopy Groups}

On the space of flagged homotopies set an equivalence relation: $\gamma
^{\prime }\sim \gamma ^{\prime \prime }$ if $\gamma ^{\prime },\gamma
^{\prime \prime }$ belong to the same connected component of $\Gamma _\sigma
^k$, i.e. there is an element $\gamma ^{k+1}\in \Gamma ^{k+1}$ such that $%
(\gamma ^{\prime },\gamma ^{\prime \prime })=(\sigma _0^k,\sigma
_1^k)(\gamma )$ (or $\gamma ^{\prime },\gamma ^{\prime \prime }$ are
cobordant). Category composition law on $\Gamma _\sigma ^k$ commutes with
this factorisation.

Take coincident highest dimension cells of $\sigma $, $\sigma
_0^{k-1}=\sigma _1^{k-1}$. Then we get a group structure on the
corresponding quotient space $\Pi _k(X,\sigma )$. 

\begin{example}
If all boundaries $(\sigma _0^i,\sigma _1^i)$ of $\sigma $ are set (belong)
to a single point $p$, then $\Pi _k(X,p)$ becomes Abelian for $k>1$, and $\Pi
_k(X,p)=\pi _k(X,p)$ is just the usual $k$-th homotopy group of manifold $X
$.
\end{example}

\subsection{Differential Calculus on $\Gamma _\sigma ^k$}

\subsubsection{Tangent to the Space of Flagged Homotopies}

For $\sigma \in \Gamma ^k\subset X$ elements of the tangent space to $%
\Gamma ^k$ at $\sigma $ are restrictions onto $\sigma $ of vector fields on $%
X$ and will be denoted $T_\sigma \Gamma ^k$.

\subsubsection{Integration over $\Gamma _\sigma ^k$}

Take a $1$-form $\varphi $ on $T_\sigma \Gamma ^k$ with values in algebra $%
{\frak g}$. For $\gamma \in \Gamma ^k$ the corresponding homotopy parameter $%
t$ is a function on $|\gamma |\subset X$. Take the vector field $v_t:=\frac{%
\partial \ln x}{\partial t}:=
\lim\limits_{\Delta t\rightarrow 0}
\frac{\ln (x_i(t+\Delta t)-x_i(t))}{\Delta t}$. 
Then we have a map $\int\exp(\varphi):\Gamma
^{k+1}\rightarrow G,\gamma \mapsto 
{\lim\limits_{\bigtriangleup t\rightarrow 0}}
\prod\limits_{i=0}^{N-1}\exp ((\varphi ,v)\Delta t)$, called {\em %
multiplicative integral}, where $G=\exp {\frak g}$, is the Lie group,
corresponding to Lie algebra of which ${\frak g}$ is the universal
envelopping one. The value $\int\limits_\gamma\exp(\varphi)$ of 
the multiplicative integration functional will also be denoted as $\varphi(\gamma)$.

\begin{example}
Take a 1-form $\omega \in \Omega ^1(X)$ with values in algebra ${\frak g}$.
It is also a $1$-form on $T\Gamma ^0$. For a point $p\subset X$ take a path $%
\gamma ([0,1])$ with $\gamma (0)=p$. Take a partition $(t_0,\ldots ,t_N)$ of
the interval $[0,1]$. Then we get a value of $\omega $ on $\gamma $ as: $%
\int\limits_\gamma\exp(\omega):=
{\lim\limits_{\bigtriangleup t\rightarrow 0}}
\prod\limits_{i=0}^{N-1}\exp ((\omega ,v(t_i))\cdot (t_{i+1}-t_i))\in G$,
where $v(t):=v(x(t))\in T_xX=T_{\gamma (t)}\Gamma ^0$ with the components 
$v_k(t):=
{\lim\limits_{\Delta t\rightarrow 0}}
\frac{\ln (x_k(t+\Delta
t)-x_k(t))}{\Delta t}$. Notice that the ordering in the ''integral product''
is fixed although the limit does not depend on a choice of parametrisation
of $\gamma $.
\end{example}

\subsubsection{Derivation in $\Gamma _\sigma ^k$}

A ${\frak g}$-valued $1$-form $\varphi $ on $T_\gamma \Gamma ^k$ is an
operator $T_\gamma \Gamma ^k\rightarrow {\frak g}$. For a pair of tangent
vectors $u,v\in T_\gamma \Gamma _\sigma ^k$ (vector fields on $\gamma $) we
can take a closed path $\gamma (u,v)$ in $\Gamma _\sigma ^k$ passing through
points $\gamma ,\gamma +\exp (u),\gamma +\exp (u)+\exp (v),{\rm and}\gamma +\exp
(u)+\exp (v)+\exp (-u)$, where $\gamma +\exp (v)$ is the result of the shift
of $\gamma $ along the integral paths of $v$ in $X$. Then take $\varphi
^{\prime }:=
{\lim\limits_{\Delta t\rightarrow 0}}\frac{\ln \varphi
(\gamma (u\Delta t^1,v\Delta t^2))}{\Delta t^1\Delta t^2}\in {\frak g}$,
where $\varphi (\gamma )$ is the multiplicative integral of $\varphi $ over $%
\gamma $. This gives us an operator $\wedge ^2T_\gamma \Gamma _\sigma
^k\rightarrow {\frak g}$ called {\em derivative} of $\varphi $.

\begin{example}
For an algebra ${\frak g}$ and $G=\exp {\frak g}$ take a $G$-valued function 
$g$ on $\Gamma ^0(X)=X$. Take a path $\gamma ([0,1])\in \Gamma ^1(X)$ with $%
\gamma (0)=p_0$ and $\gamma (1)=p_1$ and lets write $g(t):=g(\gamma (t))$.
For a partition $(t_0,\ldots ,t_N)$ of the interval $[0,1]$ we have $%
g(1)\cdot g(0)^{-1}=\prod\limits_{i=0}^{N-1}g(t_{i+1})\cdot
g(t_i)^{-1}=\prod \exp (\frac{\ln (g(t_{i+1})\cdot g(t_i)^{-1})}{\Delta t}%
\cdot \Delta t)$. If we define the directional derivative $g_\gamma ^{\prime
}:=
{\lim\limits_{\Delta t\rightarrow 0}}$ $\frac{\ln (g(t+\Delta t)\cdot
g(t)^{-1})}{\Delta t}\in {\frak g}$ then $g(1)\cdot g(0)^{-1}=
{\lim\limits_{\Delta t\rightarrow 0}}\prod \exp (g_\gamma ^{\prime }(t_i)\cdot
(t_{i+1}-t_i))$.
\end{example}

\subsection{Variational Calculus on $\Gamma _\sigma ^k$}

Take a $1$-form $\varphi \in \Omega ^1(\Gamma ^{k-1})$. It gives us a
multiplicative integration functional on $\Gamma ^k$. Take $\gamma \in
\Gamma ^k$, and for $u\in T_\gamma \Gamma _\sigma ^k$, take a variation $%
\stackrel{\sim }{\gamma }:=\gamma +\delta \gamma =\gamma +\exp (u\Delta t)$
of $\gamma $. Define the {\em first variation} of the functional as $(\delta
\varphi ,u):=
{\lim\limits_{\Delta t\rightarrow 0}}\frac{\ln (\varphi
(\gamma +\exp (u\Delta t))\cdot \varphi (\gamma )^{-1})}{\Delta t}$. It
gives us an operator $T_\gamma \Gamma _\sigma ^k\rightarrow {\frak g}$. So,
taking variational derivative of a $1$-form on $T\Gamma ^{k-1}$ we get a $1$%
-form on $T\Gamma ^k$.

\begin{proposition}
$(\delta \varphi ,u)=
{\lim\limits_{\Delta t\rightarrow 0}}\frac{\ln
((\varphi ^{\prime }\llcorner u\Delta t+\varphi )(\gamma )\cdot \varphi
(\gamma )^{-1})}{\Delta t}$.
\end{proposition}

Take a partition $(t_0,\ldots ,t_N)$ of the domain of the homotopy $\gamma $ and
the vector field $v=(\frac{\partial x}{\partial t})$. Then we can write
variation of $\varphi $ as an limit of ''integral products'' 
\begin{eqnarray*}
\delta\varphi &=&
\prod\limits_{i=0}^{N-1}\varphi(\stackrel{\sim }{\gamma(i)}
\stackrel{\sim }\gamma(i+1))\cdot 
\prod\limits_{i=0}^{N-1}\varphi(\gamma (i+1){\gamma(i)})\\
&=&
{\lim\limits_{\Delta t\rightarrow 0}}\prod\limits_{i=0}^{N-1}%
\varphi ^{\prime }\llcorner{v(x(i)+\delta
x(i))}(t_{i+1}-t_i)\cdot \prod\limits_{i=0}^{N-1}\varphi ^{\prime
}\llcorner{v(x(i))}(t_i-t_{i+1})\nonumber
\end{eqnarray*} 
Then we can expand such a product
into series in $u_i:=u(x(i))$ as $\delta \varphi =1+f_1\cdot u\Delta
t+f_2\cdot u^2\Delta t^2+o(u^2)$, where $f_k\cdot u^k$ denotes the sum of terms
of order $k$. In general this is not the expansion of 
$\exp (f_1\cdot u\Delta t)$ and since the number of commutations in terms of
order $\Delta t^2$ has order $\frac 1{\Delta t}$ then these terms make an
impact into terms of oreder $\Delta t$. Hence, unlike that in the
commutative case, the difference $\varphi (\gamma ^{\prime \prime })\cdot
\varphi (\gamma ^{\prime })^{-1}$ is not equal to $
{\lim\limits_{\Delta t\rightarrow 0}}
\prod\limits_{j=0}^{N-1}\exp ((\delta \varphi ,u)\cdot
\Delta t)$.

\subsubsection{Integrable Forms}

If the variation of the multiplicative integral functional is identically $0$
then cobordant elements of $\Gamma _\sigma ^k$ have the same value of $%
\varphi $, i.e. $\varphi (\gamma )\in G$ is constant on connected components
of $\Gamma _\sigma ^k$. Such forms will be called {\em integrable}. Since
those components correspond to elements of nonabelian homotopy group $\Pi
_k(X,\sigma )$ and the integration functional is a functor to $G$, then we have

\begin{Proposition}\label{prp:ifr}
Each integrable form gives a representation $\Pi _k(X,\sigma )\to G$. 
\end{Proposition}
The identity $\delta \varphi (u)\equiv 0$ for all $u$ implies that in the series
expansion $\delta \varphi =1+(f_1\cdot u)\Delta t+(f_2\cdot u^2)\Delta
t^2+\ldots\ $ all the coefficients $f_i=0$. In particular, $f_1\cdot u=%
{\lim\limits_{\Delta t\rightarrow 0}}\frac{\ln ((\varphi ^{\prime
}\llcorner u\Delta t+\varphi )(\gamma )\cdot \varphi (\gamma )^{-1})}{\Delta
t}\equiv 0$ implies $\varphi ^{\prime }(\gamma )\equiv 0$. Then in order to
decide on integrability of forms we have to compute their derivatives.

The value of the integral of an integrable form in general depends on the flagging $%
\sigma $.

\begin{example}
Take a connection form and an integration path in the form of a torus cut (see
Fig.3). For initial point being $A=\sigma _0=\sigma _1$ and clockwise
direction the value of the integral is 
\[
\varphi (\gamma ;A)=aba^{-1}b^{-1}
\]
For flaging set into point $B$ and the same direction the integral is 
\[
\varphi (\gamma ;B)=ba^{-1}b^{-1}a
\]

\begin{center}
\epsfxsize=8cm

Fig.3
\end{center}

\end{example}

\subsection{Integration of ${\frak g}$-valued Forms}

{\it Step 1}

Take a $1$-form $A\in \Omega ^1(X)$ with values in ${\frak g}$ (connection
form). Then it is also a $1$-form on $T\Gamma ^0$ so $A$ gives us a measure
for multiplicative integration.

\noindent
{\it Inductive step}

Suppose we have defined integration of ${\frak g}$-valued exterior $(k-1)$%
-forms on $X$. Take a sequence of forms $(\omega =\omega ^k,\ldots ,\omega
^1)$ where $\omega ^i\in \Omega ^i(X;{\frak g})$. For $\gamma \in \Gamma
^{k-1}$ and $v\in T_\gamma \Gamma _\sigma ^{k-1}$ taking a convolution of $%
\omega $ with $v$ we get a $(k-1)$-form $\omega \llcorner v+\omega ^{k-1}\in
\Omega ^{k-1}(X)$ on $\gamma \subset X$. Note, that $\omega \llcorner
v+\omega ^{k-1}$ in general is not a restriction onto $\gamma $ of any local 
$(k-1)$-form on $X$. Then we get a $1$-form $\varphi (\omega )$: $v\mapsto 
{\lim\limits_{\Delta t\rightarrow 0}}\frac{\ln ((\omega \llcorner
v\Delta t+\omega ^{k-1})(\gamma )\cdot \omega ^{k-1}(\gamma )^{-1})}{\Delta t%
}\in {\frak g}$ on $T_\gamma \Gamma _\sigma ^{k-1}$, also denoted by $\omega 
$, which is defined by induction. Since computation of $\omega ^{k-1}(\gamma
)$ is in its turn dependent of $\omega ^{k-2}$, then $\varphi $ is actually
the dependent of the total flag of forms $\varphi =\varphi (\omega ^k,\ldots
,\omega ^1)$. For $\gamma ^k\in \Gamma ^k$ take the vector field $v=(\frac{%
\partial x}{\partial t})$ corresponding to the parameter $t^k$ of homotopy.
Then we have a map $\Gamma ^k\rightarrow G,\gamma \mapsto
{\lim\limits_{\Delta t\rightarrow 0}}
\prod\limits_{j=0}^{N-1}\exp (\varphi (\omega
^k,\ldots ,\omega ^1),v\cdot \Delta t)$. Thus each element of $\Omega ^k(X;%
{\frak g})$ gives us a measure of integration on $\Gamma _\sigma ^k$. Note,
that in general the value of functional $\omega $ depends on $\gamma $. In
particular, two homotopies with the same body $|\gamma ^{\prime }|=|\gamma
^{\prime \prime }|\subset X$ but still being different paths in $\Gamma
_\sigma ^{k-1}$ may have different values of $\omega $-integral.

For each $(k-2)$-dimensional flag $\sigma $ the above defined integration
gives a functor form $\Gamma _\sigma ^k$ to $G$.

\subsection{Computation of the Exterior Derivative of a ${\frak g}$-valued Form}

We have a formal definition of exterior derivative in terms of forms on 
the space $\Gamma^k$ of paths. In this section we will see

\subsubsection{Lattice Approximation of a Homotopy}

Take a flag $\gamma =(\gamma ^1,\gamma ^2,\ldots ,\gamma ^k)$ of homotopies
and for each $\gamma ^i$ take a partition $(t_0^{(i)},\ldots ,t_{N_i}^{(i)})$
of the interval $[0,1]^{(i)}$. Each $t_j^{(i)}$ defines a hypersurface $%
\sigma (t_j^{(i)})$ in $D$, with induced flagging. The set of those
hypersurfaces for $i=0,...,k$ induces a cellular decomposition of $D$. For a
proper fineness of subdivisions those cells become homeomorphic to cubes.
This decomposition will be called {\em lattice approximation} of $D$. For
each cube $c_{n_0,\ldots ,n_k}$ take the set of its edges coming out of the
''lowest'' vertex with coordinates $(t_{n_0}^1,\ldots ,t_{n_k}^{(k)})$. An
edge with endpoint $(t_{n_0}^1,\ldots ,t_{n_i+1}^{(i)},\ldots
,t_{n_k}^{(k)})$ gives us a vector of $T_{(t_{n_0},\ldots ,t_{n_k})}D$ with
coordinates $(0,\ldots ,t_{n_i+1}^{(i)}-t_{n_i}^{(i)},\ldots ,0)$ in the
basis $\{\frac \partial {\partial t^{(0)}},\ldots ,\frac \partial {\partial
t^{(k)}}\}$. Then the frame of edges at the vertex $(n_0,\ldots ,n_k)$ gives an
element of $\wedge ^{k+1}T_{(t_{n_0},\ldots ,t_{n_k}^{(k)})}D$ - the volume
element.

\begin{proposition}
There is a cellular homomorphism $\kappa $ of a $d$-cube onto flagged $d$%
-disk.
\end{proposition}

\begin{example}
Let $d=3$. Denote by $v_{ijk},$ $e_{*ij},f_{**i}$ vertices, edges and faces
correspondingly. The homomorphism is: 

\begin{center}
$\begin{array}{rl}
v_{000}\longmapsto \sigma _0,& v_{111}\longmapsto\sigma_1\\
e_{*00}\cup e_{10*}\cup e_{1*1}\longmapsto \sigma _0^{\prime },& e_{*11}
\cup e_{01*}\cup e_{0*0}\longmapsto \sigma _1^{\prime } \\
f_{**0}\cup f_{1**}\cup f_{*1*}\longmapsto \sigma _0^{\prime \prime },& f_{**1}
\cup f_{0**}\cup f_{*0*}\longmapsto \sigma _1^{\prime \prime }
\end{array}$
\end{center}
\end{example}

This homomorphism induces a structure of flagged disk on each cubic cell $%
c_{n_0,\ldots ,n_k}$ of decomposition of $D^{k+1}$. The homotopy $\gamma
^{(k)}$ may be partitioned into a sequence of infinitesimal homotopies $%
\delta \gamma $, where each $\delta \gamma =\delta \gamma _{n_0,\ldots ,n_k}$
is a homotopy, sweeping a cubic cell $c_{n_0,\ldots ,n_k}=|\delta \gamma
_{n_0,\ldots ,n_k}|$ and concordant with its flagging structure. This
sequence induces an ordering $\Gamma $ on the set of cells. Given a value $%
\tau \in \Gamma $ we get a homotopy $\gamma _\tau ^{(k)}=\prod\limits_{\tau
^{\prime }<\tau }\delta \gamma _{\tau ^{\prime }}$ as well as the flag $%
(\gamma _\tau ^{(k)},\ldots ,\gamma _\tau ^{(0)})$ of homotopies $\gamma
_\tau ^{(i-1)}=\sigma _1(\gamma _\tau ^{(i)})\cap \sigma
(t_{n_{i-1}}^{(i-1)})$, where $\sigma _1(\gamma )$ denotes the ending border
of homotopy $\gamma $ and $\sigma (t^{(i)})$ denotes the hypersurface with
the given value $t^{(i)}$ of $i$-th coordinate on $D$.

\begin{center}
\epsfxsize=15cm

Fig.4
\end{center}

For flagged $(d+1)-$disk take its boundary $\partial D=\sigma _0^{(d)}\cup
\sigma _1^{(d)}$. The cellular homomorphism $\kappa $ of a cube onto $D$
induces a cellular subdivision of each of $\sigma _0^{(d)}$ and $\sigma
_1^{(d)}$ into $d+1$ cubes (see Fig.4). Taking pre-image of the flagging of $%
\sigma _i^{(d)}$ induces a flag of chains of faces in $\kappa ^{-1}(\sigma
_i^{(d)})$.

\begin{proposition}
There is a lattice approximation of the homotopy from the subflag $\sigma _0^{(d-1)}$
to $\sigma _1^{(d-1)}$ which may be represented as a ''connected'' path in
the set of $d-$cells of $\kappa ^{-1}(\sigma _i^{(d)})$.
\end{proposition}

Term ''connected'' here means that each next $d-$cell in the sequence has a
common $(d-1)-$cell with the previous one.

\begin{example}
Let $d=2$ and take homomorphism from the previous example. Then the
pre-image of flagged 2-sphere $\sigma _0^{\prime \prime }$ consists of $3$
faces $f_{**0},f_{1**},f_{*1*}$. $\partial \sigma _0^{\prime \prime }=\sigma
_0^{\prime }\cup \sigma _1^{\prime },\kappa ^{-1}(\sigma _0^{\prime
})=f_{*00}+f_{10*}+f_{1*1},\kappa ^{-1}(\sigma _1^{\prime
})=f_{0*0}+f_{01*}+f_{*11}$. The deformation from $\kappa ^{-1}(\sigma
_0^{\prime })$ to $\kappa ^{-1}(\sigma _1^{\prime })$ is a sequence of
chains: $f_{*00}+f_{10*}+f_{1*1},\ f_{*00}+f_{1*0}+f_{11*},\
f_{0*0}+f*_{10}+f_{11*},\ f_{0*0}+f_{01*}+f_{*11}$. We may say that each
subsequent chain in this sequence is obtained by sweeping across a $2$-face
- the formal ''difference'' between those two chains. For instance $\
f_{*00}+f_{1*0}+f_{11*}=f_{*00}+f_{10*}+f_{1*1}+f_{1**}$.
\end{example}

\subsubsection{Infinitesimal Loops in $\Gamma _\sigma ^k$}

Now take a homotopy of $\sigma _0^{(d-1)}$ to $\sigma _1^{(d-1)}$ across $%
\sigma _0^{(d)}$, and then from $\sigma _1^{(d-1)}$ to $\sigma _0^{(d-1)}$
across $\sigma _1^{(d)}$. We get a loop in the space of flagged $(d-1)$%
-disks which weeps the surface of $(d+1)$-disk $D$.

\begin{example}
For $d=3$ the $4$-disk boundary $\partial D=S^3=\sigma _0^{(3)}\cup \sigma
_1^{(3)}=D_0^3\cup D_1^3$ consists of two $3$-disks glued along a common
sphere $S^2$. In its turn $S^2=$ $\sigma _0^{(2)}\cup \sigma
_1^{(2)}=D_0^2\cup D_1^2$ glued along $S^1$. Then we can continuously drag $%
D_0^2$ to $D_1^2$ through $D_1^3$ with $S^1$ fixed and then move back to $%
D_0^2$ but now through $D_1^3$ thus sweeping the whole $\partial D$.
\end{example}

For a covering $\kappa $ of $D$ by a cube $c$ the lattice approximation of
this loop is an ordered sequence of $d$-faces of the covering cube, denoted
by $\overrightarrow{\partial }c$.

\noindent
{\it Steps 1 and 2}

For $\gamma \in \Gamma _\sigma ^k$ denote $\gamma (i):=\gamma (t_i)$. Then
for $\gamma ^{\prime },\gamma ^{\prime \prime }\in \Gamma _\sigma ^k$ and a
partition $(t_0^k,\ldots ,t_N^k)$ of the domain $[0,1]$ of homotopy we have 
\begin{eqnarray*}
\gamma ^{\prime \prime } &=&(\gamma ^{\prime }(0)\gamma ^{\prime \prime
}(0)\gamma ^{\prime \prime }(1)\gamma ^{\prime }(1)\gamma ^{\prime
}(0))\cdot (\gamma ^{\prime }(0)\gamma ^{\prime }(1))\cdot (\gamma ^{\prime
}(1)\gamma ^{\prime \prime }(1)\gamma ^{\prime \prime }(2)\gamma ^{\prime
}(2)\gamma ^{\prime }(1))\cdot (\gamma ^{\prime }(1)\gamma ^{\prime
}(2))\cdot \ldots \\
\ &=&\prod\limits_{j=0}^{N-1}(\gamma ^{\prime }(j)\gamma ^{\prime \prime
}(j)\gamma ^{\prime \prime }(j+1)\gamma ^{\prime }(j+1)\gamma ^{\prime
}(j))\cdot (\gamma ^{\prime }(j)\gamma ^{\prime }(j+1))
\end{eqnarray*}
where $(\gamma _1\gamma _2...)$ denotes the path in $\Gamma _{^{(1)}\sigma
}^{k-1}$ through a sequence of points $\{\gamma _1,\gamma _2,\ldots \}$.
Take a flag of homotopies $(\gamma ^{d+1},\gamma ^d,\ldots ,\gamma ^1)$ and
let $\gamma ^{\prime }=\gamma ^{d+1}(i),\gamma ^{\prime \prime }=\gamma
^{d+1}(i+1)$. Then we can write a lattice approximation of the above
decomposition as $\gamma (i+1)=\prod\limits_{j=0}^{N-1}\overrightarrow{%
\partial }([i,j][i+1,j+1])\cdot ([i,j][i,j+1])$, where $([i,j][i+1,j+1])$
denotes a $2$-cube (square) based on verticies $(i,j)$ and $(i+1,j+1)$, and $%
([i,j][i,j+1]):=(\gamma ^{\prime }(j)\gamma ^{\prime }(j+1))$. For a
partition $(t_0^{d-1},\ldots ,t_{N_{d-1}}^{d-1})$ of the domain of $\gamma
^{d-1}$ the $1$-loop $\overrightarrow{\partial }([i,j][i+1,j+1])$ in $\Gamma
_{^{(1)}\sigma }^{d-1}$ may be decomposed into a product: $\overrightarrow{%
\partial }([i,j][i+1,j+1])=\prod\limits_{k=0}^{N_{d-1}-1}\overrightarrow{%
\partial }([i,j,k][i+1,j+1,k+1])\cdot
([i,j,k][i,j+1,k+1])\prod\limits_{k=0}^{N_{d-1}-1}([i,j+1,k+1][i,j,k])$. The
sequence $\overrightarrow{\partial }([i,j,k][i+1,j+1,k+1])$ of cube faces is
a $2$-loop in $\Gamma _{^{(2)}\sigma }^{d-2}$.

\noindent
{\it Inductive step}

A lattice approximation of $k$-loop $\overrightarrow{\partial }%
([i^{d+1},\ldots ,i^{d-k+1}][i^{d+1}+1,\ldots ,i^{d-k+1}+1])$ in $\Gamma
_{^{(k)}\sigma }^{d-k}$ may be decomposed into a product 

$\begin{array}{l}
\prod\limits_{i^{d-k}=0}^{N_{d-k}-1}\overrightarrow{\partial }%
([i^{d+1},\ldots ,i^{d-k}][i^{d+1}+1,\ldots ,i^{d-k}+1])\cdot
([i^{d+1},i^d,\ldots ,i^{d-k}][i^{d+1},i^d+1,\ldots
,i^{d-k}+1])\\
\cdot\prod\limits_{i^{d-k}=0}^{N_{d-k}-1}([i^{d+1},i^d+1,\ldots
,i^{d-k}+1][i^{d+1},i^d,\ldots ,i^{d-k}])
\end{array}$
which is a product of $(k+1)$-paths in $\Gamma _{^{(k+1)}\sigma }^{d-(k+1)}$%
. The last term of this sequence of decompositions is 

$\begin{array}{l}
\prod\limits_{i^1=1}^{N_1-1}\overrightarrow{\partial }([i^{d+1},\ldots
,i^1][i^{d+1}+1,\ldots ,i^1+1])\cdot ([i^{d+1},i^d,\ldots
,i^1][i^{d+1},i^d+1,\ldots
,i^1+1])\\
\prod\limits_{i^1=1}^{N_1-1}([i^{d+1},i^d+1,\ldots
,i^1+1][i^{d+1},i^d,\ldots ,i^1])
\end{array}$

\noindent
which is a product of $d$-paths in $\Gamma ^0=X$.

\subsubsection{Variation of ${\frak g}$-valued Forms}

Now having a lattice approximation of the flag of homotopies $(\gamma
^{d+1},\gamma ^d,\ldots ,\gamma ^1)$ we can write $\gamma ^{d+1}(1)\cdot
\gamma ^{d+1}(0)^{-1}$ as a path through $d$-faces of $(d+1)$-cubes of
the lattice. 
For a degree $d$ flag of forms 
$\omega=(\omega^d,\dots,\omega^1)$ we have the value of
multiplicative integral functional 
\begin{eqnarray*}
&&\omega (\gamma ^{d+1}(1)\cdot \gamma ^{d+1}(0)^{-1})
=\prod_{i^{d+1}=1}^{N_{d+1}-1}(\prod_{i^d=1}^{N_d-1}\ldots \\
&&\prod_{i^1=1}^{N_1-1}\omega (\overrightarrow{\partial }([i^{d+1},\ldots
,i^1][i^{d+1}+1,\ldots ,i^1+1]))\cdot \omega ([i^{d+1},i^d,\ldots
,i^1][i^{d+1},i^d+1,\ldots ,i^1+1])\\
&&\omega ([i^{d+1},i^d+1,\ldots,i^1+1][i^{d+1},i^d,\ldots ,i^1])
\prod_{i^1=1}^{N_1-1}
\omega([i^{d+1},i^d,\ldots ,i^1][i^{d+1},i^d+1,\ldots ,i^1+1])\\
&&\ldots \prod_{i^d=1}^{N_d-1}\omega ([i^{d+1},i^d+1][i^{d+1},i^d]))
\end{eqnarray*}
where the value of $\omega$ on cells of dimension $k<d$ is the value of
the subflag $\omega^{(k)}:=(\omega^k,\dots,\omega^1)$ of corresponding dimension, 
which is defined by induction.

The integral has the following properties:
\begin{itemize}
\item
if the highest degree form $\omega^d=0$ then $\int\exp\omega=1=\exp 0\in\exp g$
\item
for a composition of integration paths $\gamma=\gamma'\circ\gamma''$ the value
of integral is $\int_\gamma\exp\omega=\int_{\gamma'}\exp\omega\cdot\int_{\gamma''}\exp\omega$
({\bf multiplicativity})
\item
if $\omega^d=d\omega^{(d-1)}$ then 
$\int_\gamma\exp\omega=\int_{\gamma(1)}\exp\omega^{(d-1)}\cdot(\int_{\gamma(0)}\exp\omega^{(d-1)})^{-1}=
\int_{d\gamma}\exp\omega^{(d-1)}$ 
({\bf Stokes formula})
\end{itemize}
If $\delta _X\omega :=\omega (\overrightarrow{\partial }([i^{d+1},\ldots
,i^1][i^{d+1}+1,\ldots ,i^1+1]))\equiv 0$ then this product contracts to $1$%
, then $\omega (\gamma ^{d+1}(1))\equiv \omega (\gamma ^{d+1}(0))$ (i.e. $%
\omega $ is integrable) and since for integrable $\omega $ this identity
does not depend on the choice of $\gamma ^{d+1}$ then the condition $\delta
_X\omega \equiv 0$ is also nesessary and $\delta _X\omega $ is called $X$-%
{\em variation} of $\omega $. Each cube $([i^{d+1},\ldots
,i^1][i^{d+1}+1,\ldots ,i^1+1])$ gives us a frame of tangent vectors $%
(v_{t_1}\Delta t_1,\ldots ,v_{t_{d+1}}\Delta t_{d+1})$ such that $%
x(i^{d+1},\ldots ,i^1)+\exp (v_{t_k}\Delta t_k)=x(i^{d+1},\ldots ,i^k,\ldots
,i^1)$. Then we have an expansion of $\delta \omega $ into series of $%
(v_{t_i})$ and the integrability $\delta \omega \equiv 0$ is equivalent to
vanishing of coefficients of all terms of this expansion. In particular, the
first (homogenious) degree terms give a map $v_{t_1},\ldots
,v_{t_{d+1}}\mapsto 
{\lim\limits_{\Delta t\rightarrow 0}}\frac{\ln
\omega (\overrightarrow{\partial }([i^{d+1},\ldots ,i^1][i^{d+1}+1,\ldots
,i^1+1]))}{\Delta t_1\ldots \Delta t_{d+1}}$ which is a ${\frak g}$-valued $%
(d+1)$-form $d\omega \in \Omega ^{d+1}(X;{\frak g})$ called {\em (first)
exterior derivative} of $\omega $. Coefficients of monomials of degree $k$
will be called {\em exterior derivatives of order} $k$.

\begin{proposition}\label{prp:dfie}
A derivative of a form is integrable form.
\end{proposition}

\begin{example}
$k=1$, $\omega =A$ is a connection form of trivial ${\frak g}$-bundle. Then
the $1$-loop is a square spanned by line elements $u,v$. The corresponding
''Riemann product'' approximation to the value of Wilson loop $A(\gamma )$
is $e^{A_u(0,0)\cdot u}\cdot e^{A_v(u,0)\cdot v}\cdot e^{-A_u(u,v)\cdot
u}\cdot e^{-A_v(0,v)\cdot v}$, where $A_u$ is the component of the
connection form in $u$ direction. Then the first derivative 
\begin{eqnarray*}
&&d\omega :=
{\lim\limits_{u,v\rightarrow 0}}\frac{\ln A(\gamma )}{u\cdot v}%
=
{\lim\limits_{u,v\rightarrow 0}}\frac{\ln (e^{A_u(0,0)\cdot u}\cdot
e^{A_v(u,0)\cdot v}\cdot e^{-A_u(u,v)\cdot u}\cdot e^{-A_v(0,v)\cdot v})}
{u\cdot v}\\
&&\ =
{\lim\limits_{u,v\rightarrow 0}}\frac{1}{u\cdot v}\ln ((1+A_uu+\frac{(A_uu)^2}{2})(1+(A_v+%
\frac{\partial A_v}{\partial u}u)v+\frac{(A_vv)^2}{2})\\
&&\ \ \ \ \ \ \cdot(1-(A_u+\frac{\partial A_u}{%
\partial u}u-\frac{\partial A_u}{\partial v}v)u+\frac{(A_uu)^2}{2})(1-(A_v-\frac{%
\partial A_v}{\partial u}u)v+\frac{(A_vv)^2}{2}))\\
&&\ =
{\lim\limits_{u,v\rightarrow 0}}\frac{\ln (1+(\frac{\partial A_u}{%
\partial v}-\frac{\partial A_v}{\partial u}+[A_u,A_v])u\cdot v)}{u\cdot v}=%
\frac{\partial A_u}{\partial v}-\frac{\partial A_v}{\partial u}+[A_u,A_v]
\end{eqnarray*}
which is the curvature form of $A$. The derivatives of higher order are
dependents of the first one, so the integrability of $1$-form is equivalent
to its being flat.
\end{example}

\begin{center}
\epsfxsize=12cm

Fig.5
\end{center}

\begin{example}
$k=2$, $(\omega ^2,\omega ^1)=(\omega ,A)$ is the $2$-flag of forms. The
lattice approximation of an infinitesimal $2$-loop is a cube spanned by 
line elements $u,v,w$. 
This loop consists of 6 moves, each sweepping a cube face (see Fig.5).
The corresponding ''Riemann product'' approximation to the integral $%
\delta \omega $ over the surface of this cube is: 
\begin{eqnarray*}
&&\omega(\overrightarrow{\partial }([0,0,0][u,v,w]))=\\
&&\omega([0,0,0][u,v,0])\cdot A([0][u])\cdot 
\omega([u,0,0][u,v,w])\cdot A([u][0])\cdot\omega([0,0,0][u,0,w])\cdot A([0][w])\\
&&\cdot \omega([0,0,w][u,v,w])\cdot A([w][0])\cdot \omega([0,0,0][0,v,w])\cdot
A([0][v])\cdot\omega([0,v,0][u,v,w])\cdot A([v][0])
\end{eqnarray*}
where each term of
the product corresponds to moving across one face of the cube in the lattice
approximation of the loop $2$-homotopy. Now, in order to find the exterior
derivatives of $\omega $ take expansions: 

$\begin{array}{l}
\begin{array}{rl}
A([u][0]):=A([u,0,0][0,0,0])=& 1-\sum_i A_iu_i-\sum_{ij}\frac{\partial A_i}{\partial x_j}u_ju_i+\frac{(\sum_iA_iu_i)^2}2+O_3(u)\\
\omega([0,0,w][u,v,w])=& 
\end{array}\\
\ \ 1+\sum_{ij}\omega_{ij}u_iv_j
+\frac{1}{2}\sum_{ijk}\frac{\partial\omega_{ij}}{\partial x_k}u_iv_jw_k
+\frac{1}{6}\sum_{ijkl}\frac{\partial^2\omega}{\partial w_k\partial w_l}u_iv_jw_kw_l
+\frac{(\sum_{ij}\omega_{ij}u_iv_j)^2}2+O_5(u,v,w)
\end{array}$\\
and so on. Then we have expansion into $u,v,w$-series of variation: 

\begin{eqnarray*}
&&\delta \omega =1+\sum_{ijk}u_iv_jw_k(\frac{\partial \omega _{ij}}{\partial x_k}+\frac{%
\partial \omega _{ki}}{\partial x_j}+\frac{\partial \omega _{jk}}{\partial x_i}%
-([\omega _{ij},A_k]+[\omega _{ki},A_j]+[\omega _{jk},A_i])) \\
&&\ \ +\sum_{ijkl}u_iu_jv_kw_l(
\frac \partial {\partial x_i}(\frac{\partial \omega _{jk}}{\partial x_l}+\frac{\partial \omega _{lj}}{\partial x_k}+\frac{\partial \omega_{kl}}{\partial x_j})
-[\omega _{kl},\frac{\partial A_i}{\partial x_j}]
+\frac{1}{2}[A_i,\frac{\partial\omega_{kl}}{\partial x_j}]\\
&&\ \ +[\omega_{il},\omega_{jk}]+\frac{\omega_{kl}A_iA_j-2A_i\omega _{kl}A_j+A_iA_j\omega_{kl}}2) \\
&&\ \ +\sum_{ijkl}u_iv_jv_kw_l(\dots)+\sum_{ijkl}u_iv_jw_kw_l(\dots)\\
&&\ \ +O_5(u,v,w)
\end{eqnarray*}
If we take $\omega =F(A)$ to be the curvature form of $A$, then the vanishing
of the coefficient at $uvw$  becomes the Bianchi identity. Note, that
vanishing of coefficients of terms of degree 4 does not follow from
vanishing of that for $(uvw)$-term, but is satisfied for curvature $F(A)$, so
these are higher analogues of Bianchi identity. 
In commutative case the coefficients of degree 4 reduce to 
$\frac \partial {\partial x_i}(\frac{\partial \omega _{jk}}{\partial x_l}+\frac{\partial \omega _{lj}}{\partial x_k}+\frac{\partial \omega_{kl}}{\partial x_j})$,
so their vanishing follows from vanishing of coefficients of degree 3.
This is why the integrability of usual differential forms may be expressed as
$d\omega=0$.
\end{example}

\subsection{Nonabelian Cohomology}

\begin{proposition}\label{prp:epr}
For integrable local section $\omega\in\Omega^{d+1}(U;{\frak g})$ of the
sheaf of $(d+1)$-forms
there is a primitive flag $\omega^{(d)}=(\omega^d,\dots,\omega^0)$
of sections $\omega^k\in\Omega^k(U;{\frak g})$ 
such that $d\omega^{(d)}=\omega^{d+1}$.
\end{proposition}

Take a local primitive $\omega_U^d\in \Omega ^d(U;{\frak g})$ of global
integrable $\omega \in \Omega ^{d+1}(X;{\frak g})$. For a path $s$, with $%
s(0)\in U$ take a sequence of open disks $\{U_i\}$ covering $s$. Then we
have a sequence of concordant local primitives $\omega _i^d\in \Omega ^d(U_i;%
{\frak g})$ of $\omega $, called {\em (analitic) continuation} of $\omega
_U^d$, which are concordant on intersections: $\omega _i^d=\omega _{i+1}^d$
on $U_i\cap U_{i+1}$. This gives a functor from $\Gamma ^1(X)$ to the
category of local sections of $\Omega ^d(X;{\frak g})$, called {\em complete 
}$d$-{\em form} (or ''multivalued form''). The space of complete $d$-forms
is denoted by $\overline{\Omega ^d(X;{\frak g})}$. If $s$ is a loop then 
we get another section $\omega _{U+s}^d\in \Omega ^d(U_i;{\frak g})$ on $U$, 
which is the result of continuation of $\omega _U^d$ along $s$. Since $%
d\omega _{U+s}^d=d\omega _U^d$ then $\delta (\omega _{U+s}^d-\omega
_U^d)\equiv 0$, and according to Proposition ~\ref{prp:epr} $\omega _{U+s}^d-\omega
_U^d=d\omega _U^{d-1}(s)$ for some $\omega _U^{d-1}\in \Omega ^{d-1}(X;%
{\frak g})$. Then for each integrable $(d+1)$-form we have a multi-valued
(complete) primitive $d$-form. The analytic continuation along a closed loop 
of the operation of taking local primitive in general sets us onto another branch
of the multivalued primitive. Vise versa the derivative of a complete $d$%
-form is an integrable univalued $(d+1)$-form, called {\em closed}. The derivatives of
single-valued forms (global sections of $\Omega ^d(X;{\frak g})$) will be
called {\em exact}. According to Proposition ~\ref{prp:ifr} we have a map 
$\overline{\Omega ^d(X;{\frak g})}\to Hom(\Pi_{d+1}(X),G)$.

{\it Question} If $\omega^i\mapsto\omega^i+\delta\omega^i$ then
what is the rule of changing highest degree form $\omega^d$ in order
for the integral to transform as 
$\int_\gamma\exp{\bar\omega}\mapsto 
\delta(\int_{\gamma^i}\exp{\bar\omega}^i)\cdot
\int_\gamma\exp{\bar\omega}\cdot
(\delta(\int_{\gamma^i}\exp{\bar\omega}^i))^{-1}$, where $\gamma^i$ is the
$0$-th border of the dimension $i$ subflag of the flagged integration 
cycle $\gamma$.

This transform gives us an action of the space $\Omega^i$ of lower degree 
form on the space $\Omega^d$ of highest order form, called {\em gauge action}.


\begin{Definition}\rm
The space of orbits of gauge action is called the {\em space of cohomologies}
with coefficients in ${\frak g}$, denoted $H^{d+1|k}(X;{\frak g})$.
\end{Definition}

\begin{example}
Take $X=S^1$,${\frak g}={\Bbb R}$, $d=0$. Then a complete $0$-form is a
multi-valued function $f$ with $d(f_{U+s}-f_U)=0$, and hence the
''discrepancy'' $\Delta :=f_{U+s}-f_U$ is constant. Adding a global section
to a comlete $0$-form preserves $\Delta $. And vise versa the difference of
two complete $0$-forms with the same discrepancy is a single-valued
function. So the space of gauge orbits of complete $0$-forms is parametrized
by values of $\Delta $, and then $H^{1|0}(S^1,{\Bbb R})={\Bbb R}^1$.
\end{example}

\begin{example}
In the framework of Section 1 take $X={\bb R}^2-\{0\}$. 
For a flat $g({\Bbb C}^*)$-valued connection form (or a locally holomorphic
function on $X$) we get the value of monodromy 
integral $\int_\gamma e^{A}\in {\Bbb C}^*$. Because of commutativity of algebra
$g({\Bbb C}^*)$ the gauge action reduces to $A\mapsto A+df$, where $f$ 
is a globally defined smooth $g({\Bbb C}^*)$-valued function. This action
preserves the monodromy, but in general moves $A$ to a different complex structure
class. Restricting this action to some fixed comlex structure we get the space
of gauge orbits 
${\cal A'}/{\cal G}=H^{1|0}({\bb R^2}-\{0\},g({\Bbb C}^*))
=\{z^\alpha|\alpha\in{\bb C}\}/\{z^n|n\in\bb Z\}={\bb C}/{\bb Z}$.
\end{example}

\section{Applications}

\subsection{De Rham theory}

Here $G={\Bbb R}_{+}^{*}$, ${\frak g=}{\Bbb R}^1$. Since $G$ is commutative
we can write all formulas in additive notations taking the logarithm of
corresponding multiplicative ones. Then in $\delta \omega $ the first
derivative is the usual exterior derivative $d\omega $ of a form and
vanishing of $d\omega $ implies vanishing of all higher derivatives. Thus
integrability is equivalent to $d\omega =0$. Since the derivative of a form
is integrable this implies the identity $d^2\omega \equiv 0$. The space of
closed $(d+1)$-forms is in 1-1 coorrespondence with complete $d$-forms
(Poicare Lemma). Then the space $H^{d+1|d}(X):=\overline{\Omega ^d(X;{\Bbb R}%
)}/\Omega ^d(X;{\Bbb R})$ of orbits of the action of $\Omega ^d(X;{\Bbb R})$
on complete $d$-forms corresponds to the De Rham cohomology space $H_{DR}^d(X;%
{\Bbb R})$. Because of commutativity the gauge action is a morphism of affine 
spaces. Then the identity $d^2\omega \equiv 0$ implies that the action of $%
\Omega ^k(X;{\Bbb R})$ on $\overline{\Omega ^d(X;{\Bbb R})}\subset
\Omega ^{d+1}(X;{\Bbb R})$ is trivial for $k<d$, so $H^{d+1|k}(X)=\overline{%
\Omega ^d(X;{\Bbb R})}$.

\subsection{Gauge theory}

Here $d=0$. The space of global $0$-forms is the space of $G$-valued
functions. It has a group structure called gauge group. Its action on the
space $\Omega ^1(X;{\frak g})$ of connections is given by the formula $%
g:A\mapsto gAg^{-1}+g^{-1}dg$. The space $\overline{\Omega ^0(X;{\frak g})}$
is in 1-1 correspondence with the space of flat connections. Since $%
d(g_{U+s}-g_U)=0$ then the discrepancy $\Delta (s):=g_{U+s}-g_U$ is
constant on $X$ similar to commutative case. A representative of 
a cohomology class of $g\in 
\overline{\Omega ^0(X;{\frak g})}$ (or of a flat connection) is given by the
set of values $\Delta (s)$ on loops $s\in \Gamma _x^1$, where $x=s(0)$, or,
since $\Delta (s)$ depends only on homotopy class of $s$, by the
representation $\pi _1(X)\rightarrow G$. Since the gauge group acts on 
$\Delta(s)$ by conjugation then $H^{1|0}(X)={\rm Hom}(\pi_1(X),G)/G$ is the space 
of conjugacy classes of representations $\{\pi_1(X)\to G\}$. The consequence $ddA\equiv 0$ of
Proposition ~\ref{prp:dfie}  is known as Bianchi identity. Note that here the action of the gauge
group does not commute with the affine space structure. In particular the action
of $\Omega ^0(X;{\frak g})$ on the space $d\Omega ^1(X;{\frak g})$ of
curvature forms is not trivial and given by $g:F\mapsto gFg^{-1}$. Then the
space of orbits $H^{2|0}(X;{\frak g}):=\overline{\Omega ^1(X;{\frak g})}%
/\Omega ^0(X;{\frak g})\neq \overline{\Omega ^1(X;{\frak g})}$ unlike that
in the commutative case.

\flushleft
\vspace{0.2 cm}
Author address:\\
College of Mathematics,\\
Independent University of Moscow

\end{document}